\documentclass[12pt]{article}
\usepackage{amsfonts}

\def\hybrid{\topmargin 0pt      \oddsidemargin 0pt
        \headheight 0pt \headsep 0pt
        \voffset=-0.5cm
        \textwidth 6.5in        
        \textheight 9in         
        \marginparwidth 0.0in
        \parskip 4pt plus 1pt   \jot = 1.5ex}
\catcode`\@=11
\def\marginnote#1{}

\newcount\hour
\newcount\minute
\newtoks\amorpm
\hour=\time\divide\hour by60
\minute=\time{\multiply\hour by60 \global\advance\minute by-\hour}
\edef\standardtime{{\ifnum\hour<12 \global\amorpm={am}%
        \else\global\amorpm={pm}\advance\hour by-12 \fi
        \ifnum\hour=0 \hour=12 \fi
        \number\hour:\ifnum\minute<10 0\fi\number\minute\the\amorpm}}
\edef\militarytime{\number\hour:\ifnum\minute<10 0\fi\number\minute}

\def\draftlabel#1{{\@bsphack\if@filesw {\let\thepage\relax
   \xdef\@gtempa{\write\@auxout{\string
      \newlabel{#1}{{\@currentlabel}{\thepage}}}}}\@gtempa
   \if@nobreak \ifvmode\nobreak\fi\fi\fi\@esphack}
        \gdef\@eqnlabel{#1}}
\def\@eqnlabel{}
\def\@vacuum{}
\def\draftmarginnote#1{\marginpar{\raggedright\scriptsize\tt#1}}
\def\draftlabel#1{{\@bsphack\if@filesw {\let\thepage\relax
   \xdef\@gtempa{\write\@auxout{\string
      \newlabel{#1}{{\@currentlabel}{\thepage}}}}}\@gtempa
   \if@nobreak \ifvmode\nobreak\fi\fi\fi\@esphack}
        \gdef\@eqnlabel{#1}}
\def\@eqnlabel{}
\def\@vacuum{}
\def\draftmarginnote#1{\marginpar{\raggedright\scriptsize\tt#1}}

\def\draft{\oddsidemargin -.5truein
        \def\@oddfoot{\sl preliminary draft {\tt(\jobname)}\hfil
        \rm\thepage\hfil\sl\today\quad\militarytime}
        \let\@evenfoot\@oddfoot \overfullrule 3pt
        \let\label=\draftlabel
        \let\marginnote=\draftmarginnote
   \def\@eqnnum{(\theequation)\rlap{\kern\marginparsep\tt\@eqnlabel}%
\global\let\@eqnlabel\@vacuum}  }

\def\numberbysection{\@addtoreset{equation}{section}
        \def\theequation{\thesection.\arabic{equation}}}

\def\titlepage{\@restonecolfalse\if@twocolumn\@restonecoltrue\onecolumn
     \else \newpage \fi \thispagestyle{empty}\c@page\z@
        \def\thefootnote{\fnsymbol{footnote}} }

\def\endtitlepage{\if@restonecol\twocolumn \else  \fi
        \def\thefootnote{\arabic{footnote}}
        \setcounter{footnote}{0}}  
\catcode`@=12
\relax

\numberbysection
\hybrid

\def\beq{\begin{equation}}
\def\eeq{\end{equation}}
\def\bea{\begin{eqnarray}}  
\def\eea{\end{eqnarray}}    
\def\p{\partial}
\def\G{\Gamma}

\def\s{\sigma}

\def\C{{\cal C}}

\def\a{\alpha}

\def\e{\varepsilon}
\def\l{\lambda}
\def\f{\varphi}
\def\s{\sigma}

\def\A{{\cal A}}
\def\B{{\cal B}}

\def\D{{\cal D}}
\def\F{{\cal F}}

\def\L{{\cal L}}

\def\O{{\cal O}}

\def\TC{{\cal T}}
\def\TD{{\cal T}^{\nu}}
\def\wt{\widetilde}
\def\wh{\widehat}

\newtheorem{theo}{Theorem}[section]

\newtheorem{cor}{Corollary}[section]
\newtheorem{lem}{Lemma}[section]

\def\bC{\mathbb C}
\def\bZ{\mathbb Z}

\def\Pic{\mathop{\rm Pic}\nolimits}

\def\res{\mathop{\rm res}\nolimits}
\def\bea{\begingroup\arraycolsep=1.5pt\begin{eqnarray}}
\def\eea{\end{eqnarray}\endgroup}

\begin{document}
 \title{Abelian solutions of the soliton equations and geometry of abelian varieties.}
\author{I. Krichever%
\thanks{Columbia University, New York, USA and
Landau Institute for Theoretical Physics, Moscow, Russia; e-mail:
krichev@math.columbia.edu. Research is supported in part by National Science
Foundation under the grant DMS-04-05519.}
\and T. Shiota%
\thanks{Kyoto University, Kyoto, Japan; e-mail: shiota@math.kyoto-u.ac.jp.
Research is supported in part by Japanese Ministry of Education, Culture, Sports,
Science and Technology under the Grant-in-Aid for Scientific Research (S)
18104001.}}

 \date{March 31, 2008}

 \maketitle

\begin{abstract} We introduce the notion of abelian solutions of the
$2D$ Toda lattice equations and the bilinear discrete Hirota equation
and show that all of them are algebro-geometric.
\end{abstract}


\section{Introduction}

The first goal of this paper to extend a theory of the abelian solutions of
the Kadomtsev-Petviashvili (KP) equation developed recently
in \cite{KrSi} to the case of the $2D$ Toda lattice
\beq\label{2DT}
\p_\xi\p_\eta\varphi_n=e^{\varphi_{n-1}-\varphi_n}-e^{\varphi_{n}-\varphi_{n+1}}
\eeq
We call a solution $\varphi_n(\xi,\eta)$
of the  equation \emph{abelian} if it is of the form
\beq\label{u}
\varphi_n(\xi,\eta)=\ln{\tau((n+1)U+z,\xi,\eta)\over
\tau(nU+z,\xi,\eta)}\,,
\eeq
where $n\in\bZ$, $\xi$, $\eta\in\bC$ and $z\in \bC^d$ are the independent variables,
$0\ne U\in\bC^d$, and for all $\xi$, $\eta$ the function
$\tau(\cdot,\xi,\eta)$ is a holomorphic section of a line bundle $\L=\L(\xi,\eta)$ on an abelian variety $X=\bC^d/\Lambda$, i.e., it satisfies the monodromy relations
\begin{equation}\label{mon0}
\tau(z+\l,\xi,\eta)=e^{a_\l\cdot z+b_\l}\tau(z,\xi,\eta), \quad \l\in\Lambda,
\end{equation}
for some $a_\l\in\bC^d$, $b_\l=b_\l(\xi,\eta)\in\bC$.

A concept of abelian solutions of soliton equations provides an unifying framework for the theory of elliptic solutions of soliton equations and the theory of their rank 1 algebro-geometric solutions. The former corresponds to the case when
the $\tau$-function is a section of line bundle on an elliptic curve ($d=1$), and the latter corresponds to the case when $X$ is the Jacobian of an auxiliary algebraic curve and $\tau$ is the corresponding Riemann $\theta$-function.

Theory of elliptic solutions of the KP equation goes back to the work \cite{amkm}, where it was found that the dynamics of poles of the elliptic
solutions of the Korteweg-de~Vries equation can be described in terms
of the elliptic  Calogero-Moser (CM) system with certain constraints.
In \cite{krelkp} it was shown that when the constraints are removed
this correspondence becomes a full isomorphism between the solutions of the elliptic  CM system and the elliptic  solutions of the KP equation.

Elliptic solutions of the $2D$ Toda lattice were considered in \cite{krzab}
where it was shown that if $\tau(z,\xi,\eta)$ in (\ref{u}) is an elliptic polynomial, i.e., if the $\tau$-function of the $2D$ Toda lattice equation is of the form
\beq\label{tauell}
\tau(z,\xi,\eta)=c(\xi,\eta)\prod_{i=1}^N \s(z-x_i(\xi,\eta))\,,
\eeq
then its zeros as functions of the variables $\xi$ and $\eta$ satisfy the equations
of motion of the Ruijsenaars-Schneider (RS) system \cite{ruij}:
$$ 
\ddot x_i = \sum_{s\neq i} \dot
x_i \dot x_s (V(x_i-x_s)-V(x_s-x_i))\,,\quad
V(x)=\zeta(x)-\zeta(x+\eta)\,,
$$ 
which is a relativistic version of the elliptic CM system.
Here and below $\s(z)=\s(z,2\omega,2\omega')$ and
$\zeta(z)=\zeta(z,2\omega,2\omega')$ are the Weierstrass $\s$- and
$\zeta$-functions, respectively.

The correspondence between finite-dimensional integrable systems and
pole systems of various soliton equations has been extensively studied
in \cite{krbab,krnest,kreltoda,krwz} (see \cite{benzvi,nekr,krlax}
and references therein for connections with the Hitchin type systems).

A general scheme of constructing Lax representations with a spectral parameter,
for systems using a specific inverse problem for
linear equations with elliptic coefficients, is presented in \cite{krnest}.
Roughly speaking, this inverse problem is the problem of characterization of
linear difference or differential equations with {\it elliptic} coefficients
having solutions that are meromorphic sections of some line bundle
on the corresponding elliptic curve (double-Bloch solutions).

Analogous problems for linear equations with coefficients that are meromorphic
functions expressed in terms of the Riemann theta function of an indecomposable
principally polarized abelian variety (ppav) $X$ were a starting point
in the recent proof in \cite{kr-schot,kr-tri} of Welters' remarkable
trisecant conjecture:
{\it an indecomposable principally polarized abelian variety $X$ is the Jacobian of a curve if and only if there exists a trisecant of its Kummer variety $K(X)$}.

Welters' conjecture, first formulated in \cite{wel1},
was motivated by Gunning's celebrated theorem \cite{gun1}
and by another famous conjecture: the Jacobians of curves are exactly
the indecomposable principally polarized abelian varieties whose theta-functions
provide explicit solutions of the KP equation. The latter
was proposed earlier by Novikov and was unsettled at the time of the Welters' work.
It was proved later in \cite{shiota}.

Let $B$ be an indecomposable symmetric matrix with positive definite imaginary part.
It defines an indecomposable principally polarized abelian variety
$X=\mathbb C^g/\Lambda$, where  the lattice $\Lambda$ is generated by the basis vectors
$e_m\in \mathbb C^g$ and the column-vectors $B_m$ of $B$.
The Riemann theta-function $\theta(z)=\theta(z|B)$ corresponding to $B$
is given by the formula
\beq\label{teta1}
\theta(z)=\sum_{m\in \mathbb{Z}^g} e^{2\pi i(z,m)+\pi i(Bm,m)},\ \
(z,m)=m_1z_1+\cdots+m_gz_g\, .
\eeq
The Kummer variety $K(X)$ is an image of the Kummer map
\beq\label{kum}
K\colon X\ni Z\longmapsto
(\Theta[\e_1,0](Z):\cdots:\Theta[\e_{2^g},0](Z))\in \mathbb{CP}^{2^g-1}
\eeq
where $\Theta[\e,0](z)=\theta[\e,0](2z|2B)$ are level two theta-functions
with half-integer characteristics $\e$.

A trisecant of the Kummer variety is a projective line which meets
$K(X)$ at least at three points. Fay's well-known trisecant formula
\cite{fay} implies that if $B$ is a matrix of $b$-periods of normalized holomorphic
differentials on a smooth genus $g$ algebraic curve $\G$, then
a set of three arbitrary distinct points on $\G$ defines a {\it one-parameter
family} of trisecants parameterized by a fourth point of the curve.
In \cite{gun1} Gunning proved under certain non-degeneracy assumptions
that the existence of such a family of trisecants
characterizes Jacobian varieties among indecomposable principally polarized
abelian varieties.

Gunning's geometric characterization of the Jacobian
locus was extended by Welters who proved that the Jacobian locus can
be characterized by the existence of a formal one-parameter family of flexes of the
Kummer varieties \cite{wel,wel1}. A flex of the Kummer variety
is a projective line which is tangent to $K(X)$ at some point up to order 2.
It is a limiting case of trisecants when the three intersection points come together.

In \cite{arb-decon} Arbarello and De Concini showed that the Welters' characterization is equivalent to an infinite system of partial differential equations  representing
the KP hierarchy, and proved that only a finite number of
these equations is sufficient. Novikov's conjecture that just the first equation of the hierarchy is sufficient for the characterization of the Jacobians is much stronger. It is equivalent to the statement that the Jacobians are characterized by the existence of length $3$ formal jet of flexes.

Welter's conjecture that requires the existence of only one trisecant
is the strongest. In fact, there are three particular cases of the Welters' conjecture, which are independent and have to be considered separately. They correspond to three possible configurations of the intersection points
$(a,b,c)$ of $K(X)$ and the trisecant:

(i) all three points coincide,

(ii) two of them coincide;

(iii) all three intersection points are distinct.

\noindent
In all of these cases the classical addition theorem for the Riemann theta-functions
directly imply that secancy conditions are equivalent to the existence
of certain solutions for the auxiliary linear problems for the KP, the $2D$ Toda, and the bilinear discrete Hirota equations, respectively.

For example, one of the Lax equations for the $2D$ Toda equation is the differential-difference equation
\beq\label{laxdd}
\p_t\psi_{n}(t)=\psi_{n+1}(t)-u_n(t)\psi_n(t)
\eeq
with the potential $u$ of the form
\beq\label{utau}
u_n(t)=\p_t\ln\tau(n,t)-\p_t\ln\tau(n+1,t)
\eeq
Let us assume that
\beq\label{tautheta}
\tau(n,t)=\theta(nU+tV+z)
\eeq
and equation (\ref{laxdd}) has a solution of the form
\beq\label{pd}
\psi_n(t)={\theta(A+nU+tV+Z)\over \theta(nU+tV+z)}\, e^{np+tE},
\eeq
where $p$, $E$ are constants and $z$ is arbitrary.
Then a direct substitution of (\ref{tautheta}) and ({\ref{pd}) into
(\ref{laxdd}) gives the equation
\beq\label{fay}
E\theta(A+z)\theta(U+z)-e^p\theta(A+U+z)\theta(z)=
\p_V\theta(U+z)\theta (A+z)-\p_V\theta(A+z)\theta(U+Z)
\eeq
which is equivalent to the condition that the projective line passing
through the points $\{K((A\pm U)/2)\}$ is tangent to the Kummer variety at
the point $ K((A-U)/2)$
(the case (ii) above).

The characterization of the Jacobian locus via (\ref{fay}) is the statement: {\it an
indecomposable, principally polarized abelian variety $(X,\theta)$ is the Jacobian of
a smooth curve of genus g if and only if there exist non-zero $g$-dimensional vectors
$U\neq A \, (\bmod\,  \Lambda),\, V$, such that equation (\ref{fay}) holds}
(\cite{kr-tri}).

The ``only if" part of the statement follows from the construction
of solutions of the $2D$ Toda lattice equations in \cite{kr-toda}, from
which it also follows that the vector $A$ in (\ref{fay}) is a point of $\G\subset J(\G)$, the vector $U$ is of the form $U=P_--P_+$, where  $P_{\pm}\in \G$ are points on $\G$, and the vector $V$ is a tangent vector
to $\G$ at one of the points.

In geometric terms the spectral curves of the elliptic RS system, that give elliptic
solutions of (\ref{2DT}) are singled out by the condition that there exist a pair of
points such that the corresponding vector $U$ {\it spans an elliptic curve} in
$J(\G)$.

For any curve $\G$ and any pair of points $P_{\pm}\in \G$
the Zariski closure of the group $\{Un|\,n\in \bZ,\
\ U=P_--P_+\}$ in $J(\G)$ is an abelian subvariety $X\subset J(\G)$. When $X$ is a
proper subvariety, i.e., $\dim X=d<g=\dim J(G)$, the restrictions of $\theta (tV+z)$
and $\theta(A+tV+z)$ on the corresponding linear subspace $\bC^d\subset \bC^g$,
i.e., the component through the origin of $\pi^{-1}(X)$, where $\pi\colon\bC^g\to J(\G)$ is the covering map, can be
seen as sections $\tau(z,t), \tau_A(z,t)$ of some line bundles on $X$, i.e. they
satisfy the monodromy properties with respect to the lattice $\Lambda\subset\bC^d$
defining $X$
\begin{equation}\label{mon}
\tau(z+\l,t)=e^{a_\l\cdot z+b_\l}\tau(z,t)\,,\ \
\tau_A(z+\l,t)=e^{a_\l\cdot z+c_\l}\tau_A(z,t)\,,\quad \l\in\Lambda,\ z\in\bC^d
\end{equation}
for some $a_\l\in\bC^d$, $b_\l=b_\l(t), c_\l=c_\l(t)\in\bC$.

Equation (\ref{fay}) restricted to $z\in \bC^d$ takes the form
\beq\label{fay1}
E\,\tau_A(z,t)\,\tau(U+z,t)-
e^p\,\tau_A(U+z,t)\,\tau(z,t)=\dot \tau(z+U,t)\,\tau_A(z,t)-
\tau(z+U,t)\,\dot\tau_A(z,t)
\eeq
Here and below ``dot" stands for the derivative with respect to the variable $t$.

At first sight equation (\ref{fay1}) considered as an equation for two unknown
sections $\tau(z,t)$ and $\tau_A(z,t)$ of some line bundles $\L(t)$ and $\L_A(t)$ on an arbitrary
abelian variety $X$ is not as restrictive as finite-dimensional equation (\ref{fay}).
Nevertheless, our first main result is that at least under certain genericity
assumptions all the {\it abelian} solutions of equation (\ref{fay1}) arise in way
described above, i.e., they are rank one algebro-geometric, and we have $X\subset
J(\G)$ for some algebraic curve $\G$, which in general might be singular.

\begin{theo} Suppose that the equation (\ref{fay1}) with some $p$, $E\in\bC$ and $0\ne U\in\bC^n$, is satisfied with $\tau(z,t),\tau_A(z,t)$, such that
for all $t$ the functions $\tau_A(z,t)$ and $\tau(z,t)$ are
holomorphic functions satisfying the monodromy properties (\ref{mon}).
Assume, moreover, that

(i) $\Lambda$ is maximal with this property, i.e., any $\l\in\bC^n$ satisfying (\ref{mon}) for some $a_\l\in\bC^n$ and $b_\l(t)$,~$c_\l(t)\in\bC$
must belong to $\Lambda$, and that,

(ii) for each $t$ the divisor $\TC^t:=\{z\in X\,\mid\tau(z,t)=0\}$ is reduced and
irreducible;

(iii) the group $\{Un|\,n\in \bZ\}$ is Zariski dense in $X$.

\noindent
Then there exist a unique irreducible algebraic curve $\G$,
smooth points $P_{\pm}\in\G$,
an injective homomorphism $j_0\colon X\to J(\G)$
and a torsion-free rank 1 sheaf
$\F\in\overline{\Pic^{g-1}}(\G)$ of degree $g-1$,
where $g=g(\G)$ is the arithmetic genus of $\G$, such that
setting $j(z)=j_0(z)\otimes\F$ we have
\beq\label{is1}
\tau(Un+z,t)=\rho(t)\,\widehat\tau_n(t,0 \mid \G,P,j(z))\,,
\eeq
where,
$\widehat \tau_n(t_1^+,t_1^-\mid \G,P,\F)$ is the $2D$ Toda
tau-function defined by the data $(\G,P_i,\F)$.
\end{theo}
Note that when $\G$ is smooth:
\beq\label{is2}
\wh\tau_n(t_1^+,t_1^-\mid\G,P,j(z))=
\theta\Bigl(nU+t_1^+V_++t_1^-V_-+j(z)\Bigm|B(\G)\Bigr)\,e^{Q(n,t_1^+,t_1^-)}\,,
\eeq
where $V_{\pm}\in\bC^n$, $Q$ is a quadratic form, $B(\G)$ is the matrix of
$B$-periods of $\G$, and $\theta$ is the Riemann theta function.
Linearization in the Jacobian $J(\G)$ of nonlinear $t$-dynamics for $\tau(z,t)$
provides some evidence that there might be underlying integrable systems on the spaces of higher level theta-functions on ppav.  The RS system is an example of
such a system for $d=1$.

Almost till the very end the proof of Theorem~1.1 goes along the lines of \cite{kr-tri}.
We would like to stress that the proof of the trisecant conjecture in \cite{kr-tri} uses
nighter of the assumptions above. We include the assumption (iii) in the statement of the
theorem only to avoid unnecessary at this stage analytical difficulties.

The second goal of this paper, discussed in the last section, is
to study abelian solutions of the BDHE. The latter is a
difference equation of the form
\beq
\tau _n (l+1,m)\tau _n (l,m+1)- \tau _n (l,m)\tau _n (l+1,m+1)+
\tau _{n+1} (l+1,m)\tau _{n-1} (l,m+1)=0\,
\label{BDHE}
\eeq
One of its auxiliary Lax equations is the two dimensional linear
difference equation
\beq\label{laxddd}
\psi(m,n+1)=\psi(m+1,n)+u(m,n)\psi(m,n)
\eeq
with the potential $u$ of the form
\beq\label{utaud}
u(m,n)={\tau(n+1,m+1)\,\tau(n,m)\over
\tau(n+1,m)\,\tau(n,m+1)}
\eeq
Under the light-cone change of variables
\beq\label{var}
x=m-n,\ \ \nu=m+n\,
\eeq
and under the assumption that $\tau(n,m)$ is of the form $\tau(Wx+z,\nu)$ with $z,W\in \bC^d$, equation (\ref{laxdd}) get transformed to the difference-functional equation
\beq\label{laxm1}
\psi(z-W,\nu)=\psi(z+W,\nu)+u\psi(z,\nu-1)\,.
\eeq
with
\beq\label{f1d5}
u(z,\nu)={\tau(z,\nu+1)\,\tau(z,\nu-1)\over
\tau(z-W,\nu)\,\tau(z+W,\nu)}
\eeq
Equation (\ref{laxm1}) for $\psi$ of the form
\beq\label{psians1}
\psi(x,\nu)={\tau_A(z,\nu)\over \tau(z,\nu)}\, e^{p\cdot z+\nu E}
\eeq
is equivalent to the discrete analog of (\ref{fay1})
\beq\label{fay2}
e^{-p\cdot W}\tau(z+W,\nu)\tau_A(z-W,\nu)=e^{p\cdot W}\tau(z-W,\nu)\tau_A(z+W,\nu)+
e^{-E}\tau(z,\nu+1)\tau_A(z,\nu-1)\,,
\eeq
where, as before, $\tau(z,\nu)$ and $\tau_A(z,\nu)$ are sections of some line
bundles on $X$, i.e. they are holomorphic functions satisfying the monodromy properties
\beq\label{monod}
\tau(z+\l,\nu)=e^{a_\l\cdot z+b_\l(\nu)}\tau(z,\nu)\,,\quad
\tau_A(z+\l,\nu)=e^{a_\l\cdot z+c_\l(\nu)}\tau_A(z,\nu)\,,\quad \l\in \Lambda,
\eeq
with respect to the lattice $\Lambda$ of an abelian variety
$X=\bC^n/\Lambda$.  If $X$ is ppav and $\tau(z,\nu)=\theta(z+V\nu),\,
\tau_A(z,\nu)=\theta(A+z+V\nu)$ then (\ref{fay2}) is equivalent to the trisecant equation
\beq\label{fay0}
e^{-p\cdot W}\theta(z+W)\theta(z+A-W)=e^{p\cdot W}\theta(z+A+W)\theta(z-W)+
e^{-E}\theta(z+V)\theta(z+A-V)\,.
\eeq
We conjecture that  under the assumption that $\tau (z,\nu),\tau_A(z,\nu)$ are
{\it meromorphic quasiperiodic functions of the variable $\nu$}
all the {\it abelian} solutions of equation (\ref{fay2})
are rank one algebro-geometric, and we have $X\subset J(\G)$ for some algebraic curve $\G$,
(which in general might be singular). The main result of the last section
is a proof of this conjecture in the case when $\tau(z,\nu)$ is periodic in the variable
$\nu$ with some sufficiently large prime period $N$. More precisely,
\begin{theo}
Suppose that the equation (\ref{fay2}) with some $p$, $E\in\bC$ and $0\ne W\in\bC^n$,
is satisfied with $\tau(z,\nu),\tau_A(z,\nu)$, such that
for all $\nu$ the functions $\tau_A(z,\nu)$ and $\tau(z,\nu)$ are
holomorphic functions satisfying the monodromy properties (\ref{monod})
with respect to the lattice $\Lambda$ of an abelian variety
$X=\bC^n/\Lambda$.  Assume, moreover, that

(i) $\Lambda$ is maximal with
this property, i.e., any $\l\in\bC^n$ satisfying (\ref{monod})
for some $a_\l\in\bC^n$ and $b_\l(\nu)$,~$c_\l(\nu)\in\bC$
must belong to $\Lambda$, and that,

(ii) for each $\nu$ the divisor $\TC^{\nu}:=\{z\in X\,\mid\tau(z,\nu)=0\}$
is reduced and is irreducible;

(iii) the Zariski closure of the group $\{2Wm|\,m\in \bZ\}$ in $X$ coincides
with $X$;

(iv) the functions $\tau(z,\nu),\tau_A(z,\nu)$ are meromorphic functions
of the variable $\nu\in \bC$ and $\tau(z,\nu)$ is a quasiperiodic
function of $\nu$, satisfying the monodromy relation
\beq\label{periodicity}
\tau(z,\nu+N)=e^{a\cdot z+c\,\nu} \tau(z,\nu)
\eeq
with an integer prime period $N>\dim H^0(\TC^{\nu})$ and with some $a\in \bC^n,\ c\in \bC$.

Then there exist a unique irreducible algebraic curve $\G$,
smooth points $P_0,P_1,P_2\in\G$,
an injective homomorphism $j_0\colon X\to J(\G)$ and a torsion-free
rank~1 sheaf $\F\in\overline{\Pic^{g-1}}(\G)$ of degree $g-1$,
where $g=g(\G)$ is the arithmetic genus of $\G$,
such that setting $j(z)=j_0(z)\otimes\F$ we have
\beq\label{is100}
\tau(Wx+z,\nu)=\rho(\nu)\,\widehat\tau(x,\nu,0,\ldots \mid \G,P,j(z))\,,
\eeq
where,
$\widehat \tau(t_1,t_2,t_3,\ldots \mid \G,P,\F)$ is the BDHE
tau-function defined by the data $(\G,P_i,\F)$.
\end{theo}

\section{Construction of the wave function}\label{sec:constr}

Equation (\ref{fay1}) is equivalent to equation (\ref{laxdd}) with
\beq\label{psitau}
u_n=-\p_t\ln{\tau((n+1)U+z,t)\over \tau(nU+z,t)},\ \
\psi_n={\tau_A(nU+z,t)\over \tau(nU+z,t)}e^{P\cdot z+Et},
\eeq
where $P\in \bC^d$ is a vector such that $P\cdot U=p$.
In the core of the proof of Theorem is the construction of quasiperiodic wave function as in (\ref{psif},\ref{psi2}) below, which contains much more information than the function $\psi$ in (\ref{psitau}) having no spectral parameter. We would like to emphasize once again that the construction of wave function follows closely the argument from the beginning of Section 2 in \cite{kr-tri} but is drastically simplified by the assumption $(ii)$  in the formulation of the theorem.

The construction is presented in two steps. First we show that the existence
of a holomorphic solutions of equation (\ref{fay2}) implies certain relations
on the tau divisor $\TC^t$.
\begin{lem} If equation (\ref{fay2}) has holomorphic solutions whose
divisors have no common components (or if the $\tau$-divisor is irreducible), then
the equation
\beq\label{rs}
\p_t^2\tau(z,t)\,\tau(z+U,t)\,\tau(z-U,t)=\p_t \tau(z,t)\,\p_t\left(\tau(z+U,t)\,\tau(z-U,t)\right)
\eeq
is valid on the divisor $\TC^t=\{\,z\in{\bC^d}\,\mid\, \tau(z,t)=0\}$.
\end{lem}
In  \cite{kr-tri} equation (\ref{rs}) was derived  with the help of
pure local consideration. Let us show that they can be easy obtained globally.

\noindent
{\it Proof.} The evaluations of (\ref{fay1}) at the divisors $\TC^t$ and $\TC^t-U$
give
\beq\label{te1}
(\dot \tau_A(z)+E\tau_A(z))\tau(z+U)=\dot \tau(z+U)\tau_A(z),\ \ z\in \TC^t,
\eeq
\beq\label{te2}
\tau_A(z)\tau(z-U)+\dot \tau(z)\tau_A(z-U)e^{-p}=0,\ \ z\in \TC^t\,.
\eeq
Here and below for brevity we omit the notations for explicit dependence of functions
on the variable $t$, i.e. $\tau(z)=\tau(z,t), \tau_A(z)=\tau_A(z,t)$.

The evaluation of  the derivative of (\ref{fay1}) at $\TC^t-U$ gives an another
equation
\beq\label{te3} (E\tau_A(z)+\dot\tau_A(z))\tau(z-U)+\dot
\tau(z-U)\tau_A(z)+\ddot\tau (z) \tau_A(z-U)e^{p}=0,\ \ z\in \TC^t
\eeq
Eliminating
$\tau_A(z-U)$ and $\dot\tau_A(z)$ from (\ref{te1}-\ref{te3}) we obtain the equation
\beq\label{te4}
\left[\ddot\tau(z)\,\tau(z+U)\,\tau(z-U)-
\dot\tau(z,t)\,\p_t\left(\tau(z+U)\,\tau(z-U)\right)\right]\tau_A(z)=0,\ \ \ z\in
\TC^t.
\eeq
which implies (\ref{rs}) due to the assumption that the divisors of $\tau$
and $\tau_A$ have no common components (or under the assumption that $\TC^t$ is irreducible).

In \cite{kr-tri} it was shown that equation (\ref{rs}) is sufficient for the existence
of {\it local} meromorphic wave solutions of (\ref{laxdd}) which are holomorphic
outside of zeros of $\tau$. Let us show that in a global setting they are sufficient
for the existence of quasi-periodic wave solutions of the differential-functional
equation:
\beq\label{laxf} \p_t\psi(z,t)=\psi(z+U,t)-u(z,t)\psi(z,t) \eeq with
\beq\label{uf} u=\p_t\ln\tau(z,t)-\p_t\ln \tau(z+U,t)\,,
\eeq
which restricted to the points
$z+Un$ takes the form (\ref{laxdd}).

The wave solution of (\ref{laxf}) is a formal solution of the form
\beq\label{psif}
\psi=k^{l\cdot z}e^{kt} \phi(z,t,k)\,,
\eeq
where $l$ is a vector $l\in \bC^d$ such that $l\cdot U=1$ and
$\phi$ is a formal series
\beq\label{psi2}
\phi(z,t,k)=e^{bt}\left(1+\sum_{s=1}^{\infty}\xi_s(z,t)\, k^{-s}\right)
\eeq

\begin{lem} Let equation (\ref{rs}) for $\tau(z,t)$ holds, and let
$\l_1,\ldots,\l_d$ be a set of linear independent vectors of the
lattice $\Lambda$ Then equation (\ref{laxf}) with $u$ as in (\ref{uf})
has a unique, up to a $z$-independent factor,
wave solution such  that:

(i) the coefficients $\xi_s(z,t)$ of the formal series (\ref{psi2})
are meromorphic functions of the variable $z\in \mathbb C^d$
with a simple pole at the divisor $\TC^t$, i.e.
\beq\label{v1}
\xi_s(z,t)={\tau_s(z,t)\over \tau(z,t)}\,;
\eeq
and $\tau_s(z,t)$ is a holomorphic function of $z$;

(ii) $\phi(z,t,k)$ is quasi-periodic with respect to the lattice $\Lambda$
\beq\label{bloch}
\phi(z+\lambda,t,k)=\phi(z,t,k)\,B^{\,\l}(k),\ \ \lambda\in \Lambda;
\eeq
and is periodic with respect to the vectors $\l_1,\ldots, \l_d$, i.e.,
\beq\label{bloch1}
B^{\,\lambda_i}(k)=1,\ \ i=1,\ldots, d.
\eeq
\end{lem}
{\it Proof.} The functions $\xi_s(z)$ are defined recursively by the equations
\beq\label{xis1}
\Delta_U\,\xi_{s+1}=\dot\xi_s+(u+b)\,\xi_s.
\eeq
Here and below
$\Delta_U$ stands for the difference derivative $e^{\p_U}-1$.
The quasi-periodicity
conditions (\ref{bloch}) for $\phi$ are equivalent to the equations
\beq\label{bloch2}
\xi_s(z+\lambda,t)-\xi_s(z,t)=\sum_{i=1}^s B^{\,\lambda}_i\xi_{s-i}(z,t)\,, \ \
\xi_0=1.
\eeq
The general quasi-periodic solution of the first equation $\Delta_U\,\xi_1=u+b$
is given by the formula
\beq\label{v5} \xi_1=-\p_t\ln \tau +l_1(z,t)\,b+c_1(t),
\eeq
where $l_1(z,t)$ is a linear form on ${\bC^d}$ such that $l_1(U,t)=1$.
It satisfies the monodromy relations (\ref{bloch2}) with
\beq\label{bloch4}
B^{\l}_1=l_1(\l)\,b -\partial_t\ln
\tau(z+\l,t)+\partial_t\ln \tau(z,t)=l_1(\l,t)\,b-\dot b_{\l}(t)\,,
\eeq
where $b_\l=b_\l(t)$ are defined in (\ref{mon}). The normalizing conditions $\B_1^{\l_i}=0,\ i=1,\ldots,d$ uniquely define the constant $b$ and
the linear form $l_1(z)$.

Let us assume that the coefficient $\xi_{s-1}$ of the series (\ref{psi2}) is known,
and that there exists a solution $\xi_s^0$ of the next equation,
which is holomorphic outside of the divisor $\TC^t$, and which
satisfies the quasi-periodicity conditions (\ref{bloch2}) with $B_s^{\l_j}=0$ and possibly $t$-dependent coefficient $B_s^\l(t)$, for $\l\neq \l_j$, i.e.
\beq\label{bloch2ur}
\xi_s(z+\lambda,t)-\xi_s(z,t)=B_s^{\l}(t)+\sum_{i=1}^{s-1} B^{\,\lambda}_i\xi_{s-i}(z,t),\ \ B_s^{\l_j}=0.
\eeq
We assume also that $\xi_s^0$
is unique up to the transformation $\xi_s=\xi_s^0+c_s(t)$, where $c_s(t)$ is
a time-dependent constant.

Let us define a function $\tau_{s+1}^0(z)$ on $\TC^t$  with the help of
the formula
\beq\label{bl1}
\tau_{s+1}^0=-\p_t \tau_s(z,t)-b\tau_s(z,t)+{\p_t\tau(z+U,t)\over\tau(z+U,t)}\tau_s(z,t), \ \
z\in \TC^t.
\eeq
Let us show that the formula (\ref{bl1}) can be written also in the alternative form:
\beq\label{bl1a}
\tau_{s+1}^0=-\p_t \tau(z,t){\tau_s(z-U,t)\over\tau(z-U,t)}, \ \
z\in \TC^t.
\eeq
By the induction assumption, $\xi_s=(\tau_s/\tau)$ is a solution of (\ref{xis1}) for
$s-1$, i.e. the function $\tau_s$ satisfies the equation
\beq\label{ind1}
[\dot \tau_{s-1}(z-U)+\tau_s(z-U)+b\tau_{s-1}(z-U)]\,\tau(z)=
\tau_s(z)\tau(z-U)+\dot\tau(z)\,\tau_{s-1}(z-U),
\eeq
where once again we omit notations for explicit dependence of all the functions on the
variable $t$.

From (\ref{ind1}) it follows that
\beq\label{ind3}
\tau_s(z)\tau(z-U)+\dot\tau(z)\,\tau_{s-1}(z-U)=0.\ \ z\in \TC^t\,.
\eeq
The evaluation of the derivative of (\ref{ind1}) at $\TC^t$ implies
\beq\label{ind2}
(\tau_s(z-U)+b\tau_{s-1}(z-U))\,\dot \tau(z)=\dot \tau_s(z)\, \tau(z-U)+\tau_s(z)\,\dot \tau(z-U)
+\ddot\tau(z)\tau_{s-1}(z-U),\ \ z\in \TC^t.
\eeq
Then, using (\ref{rs}) and (\ref{ind3}) we obtain the equation
\beq
{\dot\tau(z)\tau_s(z-U)\over \tau(z-U)}=b\tau_s(z)+\dot\tau_s(z)-{\dot \tau (z+U)\tau_s(z)\,.
\over \tau(z+U)}
\eeq
Hence, the expressions (\ref{bl1}) and (\ref{bl1a}) do coincide.

The expression (\ref{bl1}) is certainly holomorphic when
$\tau(z+U)$ is non-zero, i.e. is holomorphic
outside of $\TC^t\cap(\TC^t-U)$. Similarly from
(\ref{bl1a}) we see that $\tau_{s+1}^0(z,t)$ is
holomorphic away from  $\TC^t\cap(\TC^t+U)$.

We claim that $\tau_{s+1}^0(z,t)$ is holomorphic
everywhere on $\TC^t$. Indeed, by the assumption
the abelian subgroup generated
by $U$ is Zariski dense.
Therefore, for any point $z_0\in \TC^t$ there exists an integer $k>0$ such that $z_k=z_0-kU$
is in $\TC^t$, and $\tau(z_{k+1},t)\neq 0$. Then, from equation (\ref{bl1a}) it follows that
$\tau_{s+1}^0$ is regular at the point $z=z_k$.
Using equation (\ref{bl1}) for
$z=z_k$, we get that $\p_t\tau(z_{k-1},t)\tau_s(z_k,t)=0$. The last equality
and the equation (\ref{bl1a}) for $z=z_{k-1}$ imply that
$\tau_{s+1}^0$ is regular at the point $z_{k-1}$. Regularity of
$\tau_{s+1}^0$ at $z_{k-1}$ and
equation (\ref{bl1}) for $z=z_{k-1}$ imply $\p_t\tau(z_{k-2},t)\tau_s(z_{k-1},t)=0$.
Then equation (\ref{bl1a}) for $z=z_{k-2}$ implies that
$\tau_{s+1}^0$ is regular at the point $z_{k-2}$.
By continuing these steps
we get finally that $\tau_{s+1}^0$ is regular at $z=z_0$. Therefore, $\tau_{s+1}^0$
is regular on $\TC^t$.

Recall, that an analytic function on an analytic divisor
in $\mathbb C^d$ has a holomorphic extension onto $\mathbb C^d$ (\cite{serr}).
Therefore,
there exists a holomorphic function $\tilde\tau(z,t)$ such that
$\tilde\tau_{s+1}|_{\TC^t}=\tau_{s+1}^0$.
Consider the function $\chi_{s+1}=\tilde\tau_{s+1}/\tau$. It is holomorphic outside
of the divisor $\TC^t$.  From (\ref{bloch2}) and (\ref{bl1a}) it follows that the
function $f_{s+1}^{\l}(z)$ defined by the equation
\beq\label{bl2}
\chi_{s+1}(z+\l)-\chi_{s+1}(z)=f_{s+1}^{\l}(z)+\sum_{i=1}^s B^{\,\lambda}_i
\xi_{s+1-i}(z)\,,
\eeq
has no pole at $\TC^t$, i.e. it is a holomorphic function of $z\in \bC^d$.
It satisfies the twisted homomorphism relations
\beq\label{bl3}
f_{s+1}^{\l+\mu}(z)=f_{s+1}^{\l}(z+\mu)+f_{s+1}^{\mu}(z),
\eeq
i.e., it defines an element of the first cohomology group of $\Lambda_U$ with
coefficients in the sheaf of holomorphic functions,
$f\in H^1_{gr}(\Lambda_U,H^0(\C, \O))$. The same
arguments, as that used in the proof of the part (b) of the Lemma 12 in \cite{shiota},
show that there exists a holomorphic function $h_{s+1}(z)$
such that
\beq\label{bl4}
f_{s+1}^\l(z)=h_{s+1}(z+\l)-h_{s+1}(z)+\wt B_{s+1}^{\l},
\eeq
where $\wt B_{s+1}^{\l}=\wt B_{s+1}^{\l}(t)$ is a time-dependent constant. Hence, the function $\zeta_{s+1}=\chi_{s+1}+h_{s+1}$
has the following monodromy properties
\beq\label{bl4a}
\zeta_{s+1}(z+\l)-\zeta_{s+1}(z)=\wt B_{s+1}^{\l}+\sum_{i=1}^s B^{\,\lambda}_i
\xi_{s+1-i}(z),
\eeq
Let us consider the function
\beq\label{r2d}
R_{s+1}=\zeta_{s+1}(z+U)-\zeta_{s+1}(z)-\dot\xi_s(z)-(u(z)+b)\,\xi_s(z)
\eeq
From equation (\ref{bl1},\ref{bl1a}) it follows that it has not poles at $\TC^t$ and $\TC^t-U$, respectively. Hence, $R_{s+1}(z)$ is a holomorphic function.

From (\ref{bl4a}) it follows that it satisfies the following monodromy properties
\beq\label{rmon}
R_{s+1}(z+\l)=R_{s+1}(z)-\dot B_s^\l\,.
\eeq
Recall, that by the induction assumption $B_{s}^{\l_j}=0$, where $\l_j, j=1,\ldots,d,$ are
linear independent. Therefore, $R_{s+1}$ is a constant ($z$-independent) and
$B_s^\l$ for all $\l$ are in fact $t$-independent.

The function
\beq\label{dsol}
\wt\xi_{s+1}(z,t)=\zeta_{s+1}(z,t)+l_{s+1}(z,t)+c_{s+1}(t)\,,
\eeq
where $l_{s+1}$ is a linear form such that
$$l_{s+1}(U,t)=-R_{s+1}(t)\,,$$
is a solution of (\ref{xis1}).

Under the transformation $\xi_s\longmapsto \xi_s(z,t)+c_s(t)$ which
does not change the monodromy properties of $\xi_s$, the
solution $\wt \xi_{s+1}$ gets transformed to
\beq\label{ur2}
\xi_{s+1}=\wt \xi_{s+1}+\dot c_{s}(t)l_1(z,t)+c_s(t)\xi_1(z,t),
\eeq
where $l_1(z,t)$ is the linear form defined above in the initial step of the induction. The new solution $\xi_{s+1}$ satisfies the monodromy relations (\ref{bloch2}) with constant $B_i^\l$ for $i\leq s$ and
with $t$-dependent coefficient
\beq\label{ur1}
B_{s+1}^{\,\l}(t)=\wt B_{s+1}^{\,\l}(t)+l_{s+1}(\l,t)+\dot c_s(t)l_1(\l,t)+c_s(t)B_1^{\l}\,.
\eeq
The normalization condition (\ref{bloch1}) for $B_{s+1}^{\l_i}=1, i=0,\ldots, d$
defines uniquely $l_{s+1}$ and $\p_tc_s$, i.e. the time-dependence of $c_s(t)$.
The induction step is completed.

Note that the remaining ambiguity in the definition of $\xi_s$ on each step is the choice of a time-independent constant $c_s$. That
corresponds to the multiplication of $\psi$ by a constant formal series and thus
the lemma is proven.

\section{Commuting difference operators.}

Our next goal is to construct rings $\A^z$ of commuting difference operators
parameterized by points $z\in X$. In fact the construction of such operators
completes the proof of Theorem 1.1 because as shown in (\cite{mum,kr-dif})
there is a natural correspondence
\beq\label{corr}
\A\longleftrightarrow \{\G,P_{\pm},  \F\}
\eeq
between commutative rings $\A$ of ordinary linear
difference operators containing a pair of monic operators of co-prime orders, and
sets of algebro-geometric data $\{\G,P_{\pm}, [k^{-1}]_1,\F\}$,
where $\G$ is an algebraic curve with a fixed first jet $[k^{-1}]_1$ of a
local coordinate $k^{-1}$ in the neighborhood of a smooth
point $P_+\in\G$ and $\F$ is a torsion-free rank 1 sheaf on $\G$ such that
\beq\label{sheaf}
h^0(\G,\F(nP_+-nP_-))=h^1(\G,\F(nP_+-nP_-))=0.
\eeq
The correspondence becomes one-to-one if the rings $\A$ are considered modulo conjugation
$\A'=g(x)\A g^{-1}(x)$.

The construction of the correspondence (\ref{corr})
depends on a choice of initial point $x_0=0$. The spectral curve and the sheaf $\F$
are defined by the evaluations of the coefficients of generators of $\A$
at a finite number of points of the form $x_0+n$.
In fact, the spectral curve is independent on the choice of $x_0$, but the sheaf
does depend on it, i.e. $\F=\F_{x_0}$.

Using the shift of the initial point it is easy to show that the correspondence
(\ref{corr}) extends to the commutative rings of operators whose coefficients are
{\it meromorphic} functions of $x$. The rings of operators having poles at $x=0$
correspond to sheaves for which the condition (\ref{sheaf}) for $n=0$
is violated.

The algebraic curve $\G$ is called the spectral curve of $\A$.
The ring $\A$ is isomorphic to the ring $A(\G,P_+,P_-)$ of meromorphic functions
on $\G$ with the only pole at the points $P_+$ and which vanish at $P_-$.
The isomorphism is defined by
the equation
\beq\label{z2}
L_a\psi_0=a\psi_0, \ \ L_a\in \A, \ a\in A(\G,P_+,P_-).
\eeq
Here $\psi_0$ is a common eigenfunction of the commuting operators. At $x=0$ it is
a section of the sheaf $\F\otimes\O(P_+)$.

In order to construct rings of commutative operators we first introduce a unique
pseudo-difference operator
\beq\label{LLd}
\L(z,t)=T+\sum_{s=0}^{\infty} w_s(z,t)\,T^{-s},\ \ T=e^{\p_U},
\eeq
such that the equation
\beq\label{kkd}
\left(T+\sum_{s=0}^{N} w_s(z,t)\,T^{-s}\right)\,\psi(z,t)=k\psi(z,t)\,,
\eeq
with $\psi$ is given by (\ref{psif}), holds.
The coefficients \, $w_s(z,t)$ of $\L$ are difference polynomials in terms of the
coefficients of $\phi$. Due to quasiperiodicity of $\psi$ they are meromorphic functions on the abelian variety $X$.

Consider now the strictly positive difference parts of the operators $\L^m$.
Let $\L^m_+$ be the difference  operator such that
$\L^m_-=\L^m-\L^m_+=F_m+F_m^1T^{-1}+O(T^{-2})$. By definition the leading
coefficient $F_m$ of $\L^m_-$ is the residue of $\L^m$:
\beq\label{res1}
F_m={\rm res}_{T}\  \L^m, \ F_m^1={\rm res}_{T}\  \L^m\,T.
\eeq
From the construction of $\L$ it follows that $[\p_t-T+u, \L^n]=0$. Hence,
\beq\label{lax}
[\p_t-T+u,\L^m_+]=-[\p_t-T+u, \L^m_-]=\left(\Delta_UF_m\right)T.
\eeq
Indeed, the left hand side of (\ref{lax}) shows that the right hand side
is a difference operator with non-vanishing coefficients only at the positive
powers of $T$. The intermediate equality shows that this operator is at most
of order $1$. Therefore, it has the form $f_m T$. The coefficient $f_m$ is
easy expressed in terms of the leading coefficient $\L^m_-$. Note, that the vanishing
of the coefficient at $T^0$ and $T^{-1}$ implies the equation
\bea\label{lax55}
\Delta_U\, F_m^1&=&\p_t F_m,
\\
\Delta_U\, F_m^2&=&\p_t F_m^1+uF_1-F_1(T^{-1}u), \label{lax552}
\eea
which we will use later.

The functions $F_m(z)$ are difference polynomials in the coefficients $w_s$ of $\L$.
Hence, $F_m(z)$ are meromorphic functions on $X$.
\begin{lem} There exist holomorphic functions $q_m(z,t)$ such that
the equation
\beq\label{jq}
F_m={q_m(z+U,t)\over \tau(z+U,t)}-{q_m(z,t)\over \tau(z,t)}\,.
\eeq
holds.
\end{lem}
{\it Proof.} If $\psi$ is as in Lemma 3.1, then  there exists a unique pseudo-difference operator $\Phi$ such that
\beq\label{S}
\psi=\Phi k^{P\cdot z}e^{kt},\ \ \Phi=1+\sum_{s=1}^{\infty}\f_s(s,t)T^{-s}.
\eeq
The coefficients of $\Phi$ are universal difference polynomials in $\xi_s$.
Therefore, $\f_s(z,t)$ is a meromorphic function of $z$. Note, that $\L=\Phi\, T \Phi^{-1}$.

Consider the dual wave function defined by the left action of the operator $\Phi^{-1}$:
$\psi^+=\left(k^{-P\cdot z}e^{-kt}\right)\Phi^{-1}$.
Recall that the left action of a pseudo-difference operator is the formal adjoint action under which the left action of $T$ on a function $f$ is $(fT)=T^{-1}f$.
If $\psi$ is a formal wave solution of (\ref{laxf}),
then $\psi^+$ is a solution of the adjoint equation
\beq\label{adj}
(-\p_t-T^{-1}+u)\,\psi^+=0.
\eeq
The same arguments, as before, prove that if equation (\ref{rs})
holds then $\xi_s^+$ have simple poles on the divisor $\TC^t-U$.  Therefore, if $\psi$ as in Lemma 2.2, then the dual wave solution is of the form
$\psi^+=k^{-P\cdot z}e^{-kt}\phi^+(Ux+Z,t,k)$, where
the coefficients $\xi_s^+(z+Z,t)$ of the formal series
\beq\label{psi2+}
\phi^+(z,t,k)=e^{-bt}\left(1+\sum_{s=1}^{\infty}\xi^+_s(z,t)\, k^{-s}\right)
\eeq
have simple poles along the divisor $\TC^t-U$.

The ambiguity in the definition of $\psi$ does not affect the product
\beq\label{J0}
\psi^+\psi=\left(k^{-x}e^{-kt}\Phi^{-1}\right)\left(\Phi k^xe^{kt}\right).
\eeq
Therefore, the coefficients $J_s$ of the product
\beq\label{J}
\psi^+\psi=\phi^+(z,t,k)\,\phi(z,t,k)=1+\sum_{s=1}^{\infty}J_s(z,t)\,k^{-s}
\eeq
are meromorphic functions} on $X$. The factors in
the left hand side of (\ref{J}) have the simple poles on $\TC^t$ and $\TC^t-U$.
Hence, $J_s(z)$ is a meromorphic function on $X$ with the simple poles
at $\TC^t$ and $\TC^t-U$. Moreover,
the left and  right action of pseudo-difference operators are formally adjoint,
i.e., for any two operators the equality
$\left(k^{-x}\D_1\right)\left(\D_2k^{x}\right)=
k^{-x}\left(\D_1\D_2k^x\right)+(T-1)\left(k^{-x}\left(\D_3k^x\right)\right)$
holds. Here $\D_3$ is a pseudo-difference operator whose coefficients are difference
polynomials in the coefficients of $\D_1$ and $\D_2$. Therefore, from (\ref{J0}-\ref{201})
it follows that
\beq\label{z8}
\psi^+\psi=1+\sum_{s=1}^{\infty}J_{s}k^{-s}=
1+\Delta\left(\sum_{s=2}^{\infty}Q_sk^{-s}\right).
\eeq
The coefficients of the series $Q$ are difference polynomials in the
coefficients $\varphi_s$ of the wave operator.
Therefore, they are meromorphic functions of $z$ with poles on $\TC^t$, i.e.
$Q_s=q_s/\tau$.

From the definition of $\L$ it follows that
\beq\label{20}
\res_k\left(\psi^+(\L^n\psi)\right)k^{-1}dk=\res_k\left(\psi^+k^n\psi\right)k^{-1}dk
=J_{n}.
\eeq
On the other hand, using the identity
\beq\label{dic}
\res_k \left(k^{-x}\D_1\right)\left(\D_2k^x\right)k^{-1}dk=
\res_{T}\left(\D_2\D_1\right),
\eeq
we get
\beq\label{201}
\res_k(\psi^+\L^n\psi)k^{-1}dk=
\res_k\left(k^{-x}\Phi^{-1}\right)\left(\L^n\Phi k^x\right)k^{-1}dk=
\res_{T}\L^n=F_n.
\eeq
Therefore, $F_n=J_{n}$ and the lemma is proved.

{\bf Important remark.} In \cite{kr-tri} the statement that $F_m$
has poles only along $\TC^t$ and $\TC^t-U$
was crucial for the proof of the existence of commuting difference operators
associated with $u$. Namely, it implies that for all but a finite number
of positive integers $i\notin A$ there exist constants $c_{n,\a}$ such that
\beq\label{f1}
F_i(z,t)-\sum_{\a\in A} c_{i,\a}F_\a(z,t)=0\,,
\eeq
hence (\ref{lax}) would imply that the corresponding linear combinations
$L_i:=\L_+^i-\sum c_{i,\a}\L_+^\a$ commutes with $P:=\p_t-T-u$.
Not so: since these constants $c_{i,\a}$ might depend on $t$,
we might not have $[P,L_n]=0$, and we cannot immediately make the
next step and claim the existence of commuting operators (!).

So our next goal is to show that these constants in fact are
$t$-independent.
For that let us consider the functions $F_{i}^1(z,t)$. From (\ref{lax55}) and (\ref{jq}) it follows that
\beq\label{f1q}
F_i^1=\p_t\left({q_i(z,t)\over \tau(z,t)}\right)
\eeq
Let $\{F_\a^1\mid\a\in A\}$, for finite set $A$, be a basis of the
space $\F(t)$ spanned by $\{F_m^1\}$. Then for
all $n\notin A$ there exist constants $c_{n,\a}(t)$ such that
\beq\label{sh3}
F_n^1(z,t)=\sum_{\a\in A} c_{n,\a}(t)F_\a^1(z,t)\,.
\eeq
Due to (\ref{f1q}) it is equivalent to the equations
\bea\label{sh4}
q_n(z,t)&=&\sum_{\a} c_{n,\a}(t)q_\a(z,t)\,,\ \ z\in \TC^t,
\\
\label{sh5}
\dot q_n^1(z,t)&=&\sum_{\a} c_{n,\a}(t)\dot q_\a^1(z,t)\, \ \ z\in \TC^t\,,
\eea
from which we get
\beq\label{sh6}
\sum_{\a} (\dot c_{n,\a})q_\a(z,t)=0\,\ \ \ \ z\in \TC^t.
\eeq
From (\ref{lax552}) we obtain
\beq\label{f2eq}
\Delta_U\left(F^2_n-\sum_{\a\in A} c_{n,\a}(t)F_\a^2(z,t)\right)=
\dot c_{n,\a} F_\a^1\,.
\eeq
The left hand side is $\Delta_U$ derivative of a meromorphic function.
The right hand side has pole only at $\TC^t$. Therefore, both sides of the equation
must vanish. Then the assumption that the set $F_{\a}^1$ is minimal imply
$\dot c_{n,\a}=0$.

\begin{lem} Let $\psi$ be a wave function corresponding to $u$, and let
$L_i,\ i\notin A$  be the difference
operator given by the formula
\beq\label{a2}
L_i=\L^i_+-\sum_{\a\in A} c_{i,\a}\L^{\a}_+, \ i\notin A,
\eeq
where the constants $c_{i,\a}$ are defined by equations (\ref{sh3}).

Then the equation
\beq\label{lp}
L_i\,\psi=a_i(k)\,\psi, \ \ \ a_i(k)=k^i+\sum_{s=1}^{\infty}a_{s,i}k^{n-s}
\eeq
where $a_{s,i}$ are constants, hold.
\end{lem}
{\it Proof.} First note that from (\ref{lax}) it follows that
\beq\label{lax3}
[\p_t-T-u,L_i]=0.
\eeq
Hence, if $\psi$ is the  wave solution of (\ref{laxdd})
then $L_i\psi$ is also a wave solution
of the same equation. By uniqueness of the wave function up to a constant in $z$-factor
we get (\ref{lp}) and thus the lemma is proven.

The operator $L_i$ can be regarded as a ${z}$-parametric family   of
ordinary difference operators $L_i^{z}$.
\begin{cor} The operators $L_i^{z}$
commute with each other,
\beq\label{com1}
[L_i^{z},L_j^{z}]=0\,.
\eeq
\end{cor}
From (\ref{lp}) it follows that $[L_i^{z},L_j^{z}]\psi=0$.
The commutator is an ordinary difference operator.
Hence, the last equation implies (\ref{com1}).

\section{The fully discrete case}
The main goal of this section is to characterize under some nondegeneracy
assumptions all the abelian solutions of equation (\ref{fay2}. As above we begin
with the construction of the corresponding
quasiperiodic wave function.  We would like to emphasize once again that the construction of wave function follows closely the argument from the beginning of Section 5 in \cite{kr-tri} but is  simplified by the assumption $(iii)$  in the formulation of Theorem 1.2.

\subsection{Construction of the wave function}
First let us show that the existence of a holomorphic solutions of equation (\ref{fay2}) implies certain relations
on $\TC^\nu$.
\begin{lem} (\cite{kr-tri})\, If equation (\ref{fay2}) has holomorphic solutions, then
the equation
\beq\label{f5d}
{\tau(z+W,\nu+1)\,\tau(z-2W,\nu)\,\tau(z+W,\nu-1)\over
\tau(z-W,\nu+1)\,\tau(z+2W,\nu)\,\tau(z-W,\nu-1)}=-1\,
\eeq
is valid on the divisor ${\TC}^{\nu}=\{\,z\in{\bC^m}\,\mid\, \tau(z,\nu)=0\}$.

\end{lem}
\emph{Proof.} The evaluations of (\ref{fay2}) at the divisors $\TC^\nu\pm W$ give two different expressions for the restriction of $\tau_A(z,\nu)$ on $\TC^{\nu}$:
\beq\label{tauA1}
\tau_A(z,\nu)=e^{p\cdot W-E}\, {\tau(z+W,\nu+1)\,\tau_A(z+W,\nu-1)\over \tau(z+2W,\nu)}\,, \ \ \ z\in \TC^\nu,
\eeq
\beq\label{tauA2}
\tau_A(z,\nu)=-e^{-p\cdot W-E}\, {\tau(z-W,\nu+1)\,\tau_A(z-W,\nu-1)\over \tau(z-2W,\nu)}\,,
\ \ z\in \TC^\nu\,.
\eeq
The evaluation of equation (\ref{fay2}) for $\nu-1$ at $\TC^{\nu}$ implies
\beq\label{tauA3}
e^{-p\cdot W}\tau(z+W,\nu-1)\,\tau_A(z-W,\nu-1)=e^{p\cdot W}\tau(z-W,\nu-1)\tau_A(z+W,\nu-1),\ \ \ z\in \TC^\nu.
\eeq
Taking the ratio of (\ref{tauA1},\ref{tauA2}) and using (\ref{tauA3}) we
get  (\ref{f5d}). The lemma is proved.

Equation (\ref{f5d}) is all what we need for the rest.
\begin{lem}  Let $\tau(z,\nu)$ be a sequence of non-trivial quasiperiodic holomorphic functions on $\bC^m$. Suppose that the group $\{2W\nu|\,\nu\in \bZ\}$ is Zariski dense in $X$ and equation (\ref{f5d}) holds.
Then  there exist wave solutions $\psi(z,\nu,k)=k^\nu\phi(z,\nu,k)$
of the equation (\ref{laxm1}) with $u$ as in (\ref{f1d5})
such that:

(i) the coefficients $\xi_s(z,\nu)$ of the formal series
\beq\label{psi2d}
\phi(z,\nu,k)=\xi_0(\nu)+\sum_{s=1}^{\infty}\xi_s(z,\nu)\, k^{-s}
\eeq
are meromorphic functions of the variable $z\in \mathbb C^m$ with simple poles at the
divisor $\TD$, i.e.
\beq\label{new}
\xi_s(z,\nu)={\tau_s(z,\nu)\over \tau(z,\nu)},
\eeq
where $\tau_s(z,\nu)$ is now a holomorphic function;

(ii) $\xi_s(z,\nu)$ satisfy the following monodromy properties
\beq\label{new1}
\xi_s(z+\lambda,\nu)-\xi_{s}(z,\nu)=\sum_{i=1}^s B^{\,\lambda}_{i,\,\nu-s+i}\,\xi_{s-i}(z,\nu)\,,\ \
\ \lambda\in \Lambda,
\eeq
where $B^{\,\lambda}_{i,\,\nu}$ are $z$-independent.
\end{lem}
{\it Proof.} The functions $\xi_s(z,\nu)$ are defined recursively by the equations
\beq\label{laxm2}
\xi_{s+1}(z-W,\nu)-\xi_{s+1}(z+W,\nu)=u(z,\nu)\,\xi_s(z,\nu-1).
\eeq
The first equation for $s=-1$ is satisfied by an arbitrary $z$-independent function
$\xi_0=\xi_0(\nu)$. In what follows it will be assumed that $\xi_0(\nu)\neq 0$.

We will now prove lemma by induction in $s$. Let us assume inductively that for
$r\leq s$ the functions $\xi_r$ are known and satisfy (\ref{new1}).
Note, that the evaluation of (\ref{laxm2}) for $s-1$ and $\nu-1$ at the divisor $\TC^\nu$ gives the equation
\beq\label{feb1}
\tau_s(z-W)\tau(z+W)=\tau_s(z+W)\tau(z-W)\, ,\ \ z\in \TC^\nu.
\eeq
From (\ref{f5d}) and (\ref{feb1}) it follows that the two formulae
by which we define the residue of $\xi_{s+1}$ on $\TD$
\begin{eqnarray}
\tau_{s+1}^0(z,\nu)&=&{\tau(z+W,\nu+1)\,\tau_s(z+W,\nu-1)\over \tau(z+2W,\nu)}\,,\ \
z\in \TD,\label{may2}\\
-\tau_{s+1}^0(z,\nu)&=&{\tau(z-W,\nu+1)\,\tau_s(z-W,\nu-1)\over \tau(z-2W,\nu)}\,,\ \
z\in \TD\,. \label{may3}
\end{eqnarray}
do coincide.

The expression (\ref{may2}) is certainly holomorphic when
$\tau(z+2W)$ is non-zero, i.e. is holomorphic
outside of $\TD\cap(\TD-2W)$. Similarly from
(\ref{may3}) we see that $\tau_{s+1}^0(z,\nu)$ is
holomorphic away from  $\TD\cap(\TD+2W)$.

We claim that $\tau_{s+1}^0(z,\nu)$ is holomorphic
everywhere on $\TD$. Indeed, by assumption
the closure of the abelian subgroup generated
by $2W$ is everywhere dense. Thus for any $z\in\TD$
there must exist some $N\in\mathbb N$ such that
$z-2(N+1)W\not\in\TD$; let $N$ moreover be the
minimal such $N$. From (\ref{may3}) it then follows that
$\tau_{s+1}^0(z,\nu)$ can be extended holomorphically to the point
$z-2NW$. Thus expression
(\ref{may2}) must also be holomorphic at $z-2NW$; since its
denominator there vanishes, it means that the numerator must also
vanish.
But this expression is equal to the numerator of (\ref{may3}) at
$z-2(N-1)W$; thus $\tau_{s+1}^0$ defined from
(\ref{may3}) is also holomorphic at $z-2(N-1)W$ (the
numerator vanishes, and the vanishing order of the denominator is
one, since we are talking exactly about points on its vanishing
divisor).
Note that we did not quite need the fact $z-2(N+1)W\not\in\TD$
itself, but the consequences of the minimality of $N$, i.e.,
$z-2kW\in\TD$, $0\leq k\leq N$,
and the holomorphicity of $\tau_{s+1}^0(z,\nu)$ at $z-2NW$."
Therefore, in the same
way, by replacing $N$ by $N-1$, we can then deduce holomorphicity
$\tau_{s+1}^0(z,\nu)$ at $z-2(N-2)W$ and, repeating the process $N$ times, at
$z$.

Recall that an analytic function on an analytic divisor
in $\mathbb C^d$ has a holomorphic extension to all of $\mathbb C^d$ (\cite{serr}).
Therefore, there exists a holomorphic function $\widetilde \tau_{s+1}(z,\nu)$
extending the $\tau_{s+1}^0(z,\nu)$. Consider then the function
$\chi_{s+1}(z,\nu)=\wt\tau_{s+1}(z,\nu)/\tau (z,\nu)$,
holomorphic outside of $\TD$ .

From (\ref{new1}) and (\ref{may2}) it follows that the function
\beq\label{new6}
f_{s+1}^{\l}(z,\nu)=\chi_{s+1}(z+\l,\nu)-
\chi_{s+1}(z,\nu)-\sum_{i=1}^{s} B^{\,\lambda}_{i,\,\nu-1-s+i}\,
\xi_{s+1-i}(z,\nu)
\eeq
has no pole at the divisor $\TD$. Hence, it  is a holomorphic function.
It satisfies the twisted homomorphism relations
\beq\label{new3}
f_{s+1}^{\l+\mu}(z,\nu)=f_{s+1}^{\l}(z+\mu,\nu)+f_{s+1}^{\mu}(z,\nu),
\eeq
i.e., it defines an element of the first cohomology group of $\Lambda_0$ with
coefficients in the sheaf of holomorphic functions,
$f\in H^1_{gr}(\Lambda_0,H^0(\mathbb C^m, \O))$. Once again using the same
arguments, as that used in the proof of the part (b) of the Lemma 12 in \cite{shiota},
we get that there exists a holomorphic function $h_{s+1}(z,\nu)$
such that
\beq\label{new4}
f_{s+1}^{\l}(z,\nu)=h_{s+1}(z+\l,\nu)-h_{s+1}(z,\nu)+\wt B_{s+1,\,\nu}^{\l}\xi_0(\nu),
\eeq
where $\wt B_{s+1,\,\nu}^{\l,}$ is $z$-independent.
Hence, the function $\zeta_{s+1}=\chi_{s+1}+h_{s+1}$
has the following monodromy properties
\beq\label{new5}
\zeta_{s+1}(z+\l,\nu)-\zeta_{s+1}(z,\nu)=
\wt B_{s+1,\nu}^{\l}\,\xi_0(\nu)
+\sum_{i=1}^{s} B^{\,\lambda}_{i,\,\nu-1-s+i}\,\xi_{s+1-i}(z,\nu)\,.
\eeq
Let us consider the function $R_{s+1}$
defined by the equation
\beq\label{R}
R_{s+1}=\zeta_{s+1}(z-W,\nu)-\zeta_{s+1}(z+W,\nu)-
u(z,\nu)\,\xi_{s}(z,\nu-1)\,.
\eeq
Equation (\ref{may2}) and (\ref{may3}) imply that the r.h.s of (\ref{R}) has
no pole at $\TD\pm W$. Hence,
$R_{s+1}(z,\nu)$ is a holomorphic function of $z$.
From (\ref{new1},\ref{new5}) it follows that it is periodic
with respect to the lattice $\Lambda$, i.e., it is a function on $X$.
Therefore, $R_{s+1}$ is a constant.

Hence, the function
\beq\label{sol}
\xi_{s+1}(z,\nu)=\zeta_{s+1}(z,\nu)+l_{s+1}(z,\nu)\xi_0(\nu)+c_{s+1}(\nu)\xi_0(\nu)\,,
\eeq
where $c_{s+1}(\nu)$ is a constant, and $l_{s+1}$ is a linear form such that
$$l_{s+1}(2W,\nu)\xi_0(\nu)=-R_{s+1}(\nu)\,,$$
is a solution of (\ref{laxm2}). It satisfies the monodromy relations (\ref{new1})
with
\beq\label{new7}
B_{s+1,\,\nu}^{\,\l}=\wt B_{s+1,\,\nu}^{\,\l}+l_{s+1}(\l,\nu)\,.
\eeq
The induction step is completed and thus the lemma is proven.

On each step the ambiguity in the construction of $\xi_{s+1}$ is
a choice of linear form $l_{s+1}(z,\nu)$ and constants $c_{s+1}(\nu)$. As
in Section 2, we would like to fix this ambiguity by normalizing monodromy coefficients $B_{i,\,\nu}^{\l}$ for a set of linear independent vectors
$\l_1,\ldots,\l_d\in\Lambda$. As it was revealed in \cite{kr-tri} in the fully discrete case there is an obstruction for that.
This obstruction is a possibility of the existence of periodic solutions
of (\ref{laxm2}), $\xi_{s+1}(z+\l,\nu)=\xi_{s+1}(z,\nu), \ \l\in \Lambda,$
for $0\leq s\leq r-1$.

Note, that there are no periodic solutions of (\ref{laxm2}) for all $s$. Indeed,
the functions $\xi_s(z,\nu)$ as solutions of non-homogeneous equations are linear
independent. Suppose not.  Take a smallest nontrivial linear relation among
$\xi_s(z,\nu)$, and apply (5.24) to obtain a smaller linear relation.
The space of meromorphic functions on $X$
with simple pole at $\TD$ is finite-dimensional. Hence, there exists minimal
$r$ such that equation (\ref{laxm2}) for $s=r$ has no periodic solutions.

Let $\l_1,\ldots,\l_d$ be a set of linear independent vectors
in $\Lambda$. Without loss of generality throughout the rest of the paper it
will be assumed that there is no linear form $l(z),\ z\in \bC^m,$  with $l(\l_j)=1$ and $l(2W)=0$.

\begin{lem} Suppose equations (\ref{laxm2}) has periodic solutions for $s<r$
and has a quasi-periodic solution $\xi_{r}$ whose monodromy relations for
$\l_j$ have the form
\beq\label{new10}
\xi_{r}(z+\l_j,\nu)-\xi_{r}(z,\nu)=b\,\xi_0(\nu),\ \ j=1,\ldots,d,
\eeq
where $b\neq 0$ is a constant.
Then for all $s$ equations (\ref{laxm2}) has solutions of the form (\ref{new})
satisfying (\ref{new1}) with $B_{i,\,\nu}^{\l_j}=b\,\delta_{i,r}$,
i.e.,
\beq\label{new20}
\xi_{s}(z+\l_j,\nu)-\xi_{s}(z,\nu)=b\,\xi_{s-r}(z,\nu).
\eeq
\end{lem}
{\it Proof.} We will now prove the lemma by induction in $s\geq r$. Let us assume
inductively that $\xi_{s-r}$ is known, and for $1\leq i\leq r$ there are solutions
$\tilde \xi_{s-r+i}$ of (\ref{laxm2}) satisfying
(\ref{new1}) with $B_{i,\,\nu}^{\l_j}=b\,\delta_{i,r}$. Then,
according to the previous lemma, there exists a solution $\tilde\xi_{s+1}$
of (\ref{laxm2}) having the form (\ref{new}) and  satisfying monodromy relations
(\ref{new1}), which for $\l_j$ have the form
\beq\label{new30}
\tilde \xi_{s+1}(z+\lambda_j,\nu)-\tilde \xi_{s+1}(z,\nu)=
b\,\tilde \xi_{s-r+1}(z,\nu)+
B_{s+1,\,\nu}^{\l_j}\xi_0(\nu)\,.
\eeq
If $\xi_{s-r}$ is fixed, then the general quasi-periodic solution $\xi_{s-r+1}$
with the normalized monodromy relations is of the form
\beq\label{transc}
\xi_{s-r+1}(z,\nu)=\wt \xi_{s-r+1}(z,\nu)+c_{s-r+1}(\nu)\xi_0(\nu)\,.
\eeq
It is easy to see that under the transformation (\ref{transc}) the functions
$\wt \xi_{s-r+i}$ get transformed to
\beq\label{new40}
\xi_{s-r+i}(z,\nu)=
\wt \xi_{s-r+i}(z,\nu)+c_{s-r+1}(\nu-i+1)\, \xi_{i-1}(z,\nu)\,.
\eeq
This transformation does not change the monodromy properties of $\xi_{s-r+i}$
for $i\leq r$, but changes the monodromy property of $\xi_{s+1}$:
\bea
\xi_{s+1}(z+\lambda_j,\nu)-\xi_{s+1}(z,\nu)&=&b\,\xi_{s-r+1}(z,\nu)+
B_{s+1,\,\nu}^{\l_j}\,\xi_0(\nu)+ \nonumber
\\
&{}&b\,(c_{s-r+1}(\nu-r)-c_{s-r+1}(\nu))\,\xi_0(\nu) \label{new31}.
\eea
Recall, that $\wt \xi_{s+1}$ was defined up to a linear form $l_{s+1}(z,\nu)$
which vanishes on $2W$. Therefore the normalization of
the monodromy relations for $\xi_{s+1}$ uniquely defines this form and
the differences $(c_{s-r+1}(\nu-r)-c_{s-r+1}(\nu))$. The induction
step is completed and the lemma is thus proven.

Note, the following important fact: if $\xi_{s-r}$ is fixed
then $\xi_{s-r+1}$, such that there exists quasi-periodic solution $\xi_{s+1}$
with normalized monodromy properties,
is defined uniquely up to the transformation:
\beq\label{transcc}
\xi_{s-r+1}(z,\nu)\longmapsto\xi_{s-r+1}(z,\nu)+c_{s-r+1}(\nu)\xi_0(\nu),
\ \ c_{s-r+1}(\nu+r)=c_{s-r+1}(\nu).
\eeq
Our next goal is to show that the assumption of Lemma 4.3 holds for some $r$, and
then to fix the remaining ambiguity (\ref{transcc}) in the definition of
the wave function. At this moment we are going to use for the first time
the assumption that $\tau$ is a meromorphic periodic function of the variable $\nu$.

Let $r$ be the minimal integer such that there exist solutions
$\xi_0^0=1,\xi_1^0,\ldots,\xi_{r-1}^0$ of (\ref{laxm2})
that are periodic functions of $z$ with respect to
$\Lambda$, and there is no periodic solution $\xi_r$ of (\ref{laxm2}).
As it was noted above, the functions $\tau_s$ are linear independent. Hence,
$r\leq h^0(Y,\theta|_Y)$.

If $\xi_{r-1}^0$ is periodic, then the monodromy relation for
$\xi_r$ has the form
\beq\label{new65}
\xi_r^0(z+\lambda,\nu)-\xi_{r}^0(z,\nu)=B^{\,\lambda}_{r}(z,\nu)\,,\ \
\ \lambda\in \Lambda.
\eeq
The function $B_r^\l$ is independent of the ambiguities in the definition
of $\xi_i, \ i<r$, and therefore, it
is a well-defined holomorphic function of $z\in X$.
Hence, it is $z$-independent, $B_{r}^\l(z,\nu)=B_r^\l(\nu)$.
The function $\xi_r^0$ is defined up to addition of a linear form $l_r(z,\nu)$
such that $l(2W,\nu)=0$. Therefore, there exist the solution $\xi_r^0$
such that $B_r^{\l_j}(\nu)=B_r(\nu)$. There is no $\xi_r^0$ which is periodic
for all $\nu$. Hence, $B_r(\nu)\neq 0$ at least for
one value of $\nu$.  By assumption the function $\tau$ is a meromorphic function of $\nu$. Therefore, $B_r(\nu)$ is a meromorphic function of $\nu$. Shifting $\nu\to \nu+\nu_0$ if needed, we may assume without loss of generality that
$B_r(\nu)\neq 0$ for {\it all} $\nu\in \bZ$. From (\ref{periodicity}) it follows that
$u(z,\nu+N)=u(z,\nu)$. Hence, $B_r(\nu)$ is a periodic function of $\nu$, i.e.
\beq\label{perB}
B_r(\nu+N)=B_r(\nu).
\eeq
Under the transformation
\beq\label{ttran0}
\xi_0^0=1\longmapsto \xi_0(\nu)
\eeq
the solutions $\xi_r^0$ get transformed to
\beq\label{ttrans}
\xi_s(z,\nu)=\xi_s^0(z,\nu)\,\xi_0(\nu-s).
\eeq
From (\ref{new65}) it follows that the transformed function $\xi_r$ satisfies the
relations
\beq\label{newtransf}
\xi_r(z+\lambda,\nu)-\xi_{r}(z,\nu)=B^{\,\lambda}_{r}(\nu)\xi_0(z,\nu-r)\,,\ \
\ \lambda\in \Lambda.
\eeq
The equation
\beq\label{eqforxi0}
b\,\xi_0(\nu)=B_r(\nu)\xi_0(\nu-r), \ \ \xi_0(\nu+N)=\xi_0(\nu).
\eeq
restricted to the space of periodic functions $\xi_0$ can be regarded as a
finite-dimensional linear equation. The vanishing of the determinant of this equation
defines the constant $b$. With $b$ fixed equation (\ref{eqforxi0}) defines
$\xi_0$ uniquely up to multiplication by a function $c_0(\nu)$ such that
$c_0(\nu+N)=c_0(\nu+r)=c_0(\nu)$. By the assumption of Theorem 1.2 the period $N$
is prime and $N>H^0(\TC^\nu)$. As it was mentioned above $r\leq H^0(\TC^\nu)$.
Hence, two periods of $c_0$ are coprime, i.e.,$(r,N)=1$.
Therefore, $\xi_0$ is defined uniquely up to a constant factor.

\begin{lem} Suppose that the assumptions of Theorem 1.2 hold.
Then  there exists a formal solution
\beq\label{ff1}
\phi=\xi_0(\nu)+\sum_{s=1}^{\infty}\xi_s(z,\nu)\,k^{-s}
\eeq
of the equation
\beq\label{ff2}
k\phi(z-W,\nu,k)=k\phi(z+W,\nu,k)+u(z,\nu)\,\phi(z,\nu-1,k)\,,
\eeq
with $u$ as in (\ref{f1d5}) such that:

(i) the coefficients $\xi_s$ of the formal series $\phi$
are of the form $\xi_s=\tau_s/\theta$, where
$\tau_s(Z)$ are holomorphic functions;

(ii) $\phi(z,\nu, k)$ is quasi-periodic with respect to the lattice $\Lambda$ and
for the basis vectors $\l_j$ in $\bC^m$ its monodromy relations
have the form
\beq\label{ff4}
\phi(z+\lambda_j,\nu,k)=(1+b\,k^{-r})\,\phi(z,\nu,k),\ \ j=1,\ldots, m,
\eeq
where $b$ are constants defined by (\ref{eqforxi0});

(iii) $\phi(z,\nu, k)$ is a quasi-periodic function of the variable $\nu$, i.e.
\beq\label{pernu}
 \phi(z,\nu+N, k)=\phi(z,\nu, k) \mu(k)
 \eeq
(iv) $\phi$ is unique up to the multiplication by a constant in $z$ factor
$\rho(k)$.
\end{lem}
{\it Proof.} We prove the lemma by induction in $s$. Let us assume
inductively that $\xi_{s-r}$ is known. As shown above the normalization of the relations for $\xi_{s+1}$ uniquely defines $\xi_{s-r+1}$ up to the transformation
(\ref{transcc}), i.e. up to a $r$-periodic function $c_{s-r+1}(\nu+r)=c_{s-r+1}(\nu)$.
The quasiperiodicity condition (iii)  is equivalent to the
condition that this function of $c_{s-r+1}$ is $N$-periodic. As it was mentioned above the periods $r$ and $N$ are coprime. Hence, on each step $\xi_{s-r+1}$
is defined up to the additive constant. This ambiguity corresponds to the
multiplication of $\phi$ be a constant factor $\rho(k)$, and thus the lemma
is proven.

\subsection{Commuting difference operators}

As it Section 3 we are now going to construct
rings $\A^z$ of commuting difference operators.
First we introduce pseudo-difference operator in one of the original variable $m$ depending on the second variable $n$
and a point $z\in \bC^d$. (Recall, that the variables $n,m$ are related to $x,\nu$ via (\ref{var}).

The formal series $\phi(z,\nu,k)$ defines a unique
pseudo-difference operator
\beq\label{LLdd}
\L(z,\nu)=w_0(\nu)T+\sum_{s=0}^{\infty} w_{s+1}(z,\nu)\,T^{-s},\ \ T=e^{\p_m},
\eeq
such that the equation
\beq\label{kkdd}
\left(w_0(m+n)T+\sum_{s=0}^{N} w_s(z+(m-n)W,(m+n))\,T^{-s}\right)\,\psi=k\psi\,.
\eeq
holds. Here $\psi=k^{n+m}\phi(z+(m-n)W,(m+n),k)$.
The coefficients \, $w_s(z,\nu)$ of $\L$ are difference polynomials in terms of the
coefficients of $\phi$. Due to quasiperiodicity of $\psi$ they are meromorphic functions on the abelian variety $X$.

From equations (\ref{ff2},\,\ref{kkdd}) it follows that
\beq\label{m22}
\left((\Delta_1\L^i)\,T_1-(\Delta \L^i) \,T-[u,\L^i]\right)\psi=0\,,
\eeq
where $\Delta_1\L^i$ and $\Delta\L^i$ are pseudo-difference operator in $T$, whose
coefficients are difference derivatives of the coefficients of $\L^i$
in the variables $n$ and $m$ respectively.
Using the equation $(T_1-T-u)\,\psi=0$, we get
\beq\label{m23}
\left(\left(\Delta_1\L^i\right)T-\left(\Delta\L^i\right) T +
\left(\Delta_1 \L^i\right)u-[u,\L^i]\right)\psi=0.
\eeq
The operator in the left hand side of (\ref{m23}) is a pseudo-difference operator
in the variable $m$. Therefore, it has to be equal to zero. Hence, we have the equation
\beq\label{m21}
\left(\Delta_0\L^i\right) T +\left(\Delta_1 \L^i\right)u-[u,\L^i]=0, \
\Delta_{0}=T_1-T
\eeq
Let $\L^i_+$ be the strictly positive difference part of the operator $\L^i$, i.e.,
\beq\label{l+}
\L^i=\L_+^i+\L^i_-=\L_+^i+\sum_{s=0}^{\infty}F_{i,s}T^{-s}
\eeq
Then,
\beq\label{m24}
\left(\Delta_{0}\L^i_+\right) T +\left(\Delta_1 \L^i_+\right)u-[u,\L^i_+]=
-\left(\Delta_{0}\L^i_-\right) T -\left(\Delta_1 \L^i_-\right)u+[u,\L^i_-]
\eeq
The left hand side of (\ref{m24})
is a difference operator with non-vanishing coefficients only at the positive
powers of $T$. The right hand side is a pseudo-difference operator
of order $1$. Therefore, it has the form $f_i T$. The coefficient $f_i$ is
easy expressed in terms of the leading coefficient $\L^i_-$. Finally we get
the equation
\beq\label{m25}
\left(\Delta_{0}\L^i_+\right) T +\left(\Delta_1 \L^i_+\right)u-[u,\L^i_+]=
-(\Delta_{0}F_i)\,T,
\eeq
where $F_i=F_{i}=\res\  \L^i$.

By definition of $\L$ we have that the functions $F_{i}$ in (\ref{l+}) are of the form
\beq\label{l+1}
F_{i}=\res_T \L^i=F_{i}(z+(m-n)W,(m+n))
\eeq
where for each $\nu$ the functions $F_{i}(z,\nu)$ are abelian
functions, i.e., periodic functions of the variable $z\in \bC^d$.

\begin{lem} The abelian functions $F_i$ have the form
\beq\label{fiq}
F_i(z,\nu)={q_i(z+W,\nu+1)\over \tau(z+W,\nu+1)}-
{q_i(z,\nu)\over \tau(z,\nu)}\,,
\eeq
where $q_i(z,\nu)$ are holomorphic functions of the variable $z\in \bC^d$.
\end{lem}
\noindent
{\it Proof.}
The wave solution $\psi$ define the unique operator $\Phi$ such that
\beq\label{Sd}
\psi=\Phi k^{n+m},\ \
\Phi=1+\sum_{s=1}^{\infty}\f_s((m-n)W+z,(m+n)\,T^{-s}\,,
\eeq
where $\f_s(z,\nu)$ are meromorphic functions of $z\in \bC^d$.
The dual wave function
\beq\label{psinewd}
\psi^+=k^{-n-m}\left(1+\sum_{s=1}^{\infty}\xi^+_s((n-m)W+z,(n+m))\, k^{-s}\right)
\eeq
is defined by the formula
\beq\label{m26}
\psi^+=k^{-n-m}\, T_1\,\Phi^{-1}\,T_1^{-1}.
\eeq
It satisfies the equation
\beq\label{adjd}
(T_1^{-1}-T^{-1}-u)\,\psi^+=0,
\eeq
which implies that the functions $\xi_s^+(z,\nu)$ have
the form $\xi_s^+(z,\nu)=\tau_s^+(z,\nu)/\theta(z+W,\nu+1)$,
where $\tau_s^+(z,\nu)$ are holomorphic functions of $z\in \bC^d$.
Therefore, the functions $J_s(z,\nu)$ such that
\beq\label{Jnew}
(\psi^+T_1)\,\psi=k+\sum_{s=1}^{\infty}J_s((n-m)W+z,(n+m))\,k^{-s+1}
\eeq
are meromorphic function on $X$ with the simple poles at $\TC^{\nu}$ and $\TC^{\nu+1}-W$.

The same arguments as that used for the proof of (\ref{z8}) show that
\beq\label{z8d}
(\psi^+T_1)\psi=(k^{-x}T_1\Phi^{-1})(\Phi k^x)=k+(\Delta Q)
\eeq
where the coefficients of the series $Q$ are of the form
\beq\label{Q}
Q=\sum_{s=0}^{\infty} Q_s(n-m)W+z,(n+m)) k^{-s}\,,
\eeq
and the functions $Q_s(z,\nu)$ are difference polynomials in the
coefficients $\varphi_s$ of the wave operator. Therefore, they are well-defined
meromorphic functions of $z$.
As shown above, the functions
\beq\label{JQ}
J_s(z,\nu)=Q_s(z+W,\nu+1)-Q_s(z,\nu)
\eeq
have simple poles at $\TC^{\nu}$ and $\TC^{\nu+1}-W$. Hence,
$Q_s(z,\nu)$ have poles only at $\TC^\nu$, i.e.
\beq\label{QQ}
Q_s={q_s(z,\nu)\over \tau(z,\nu)}\,,
\eeq
where $q_s(z,\nu)$ are holomoprhic functions of $z$.

From the definition of $\L$ it follows that
\beq\label{20new}
\res_k\left((\psi^+T_1)\,(\L^i\psi)\right)k^{-2}dk=
\res_k\left((\psi^+\,T_1)\,\psi\right)k^{i-2}dk
=J_{i}.
\eeq
On the other hand, using (\ref{dic})
we get
\beq\label{20nn}
\res_k((\psi^+\,T_1)\,(\L^i\psi)\,k^{-2}dk=
\res_k\left(k^{-n-m}\Phi^{-1}\right)\left(\L^i\Phi k^{n+m}\right)k^{-1}dk=
\res_{T}\L^i=F_i.
\eeq
Equation (\ref{fiq}) is a direct corollary of (\ref{JQ}-\ref{20nn}).
The lemma is proved.

The function $\psi$ is quasiperiodic function of the variable $\nu$. Then, from the definition
of $\psi^+$ it follows that
\beq\label{pernu+}
 \phi^+(z,\nu+N, k)=\phi^+(z,\nu, k) \mu^{-1}(k)\,,
 \eeq
where $\mu(k)$ is defined in (\ref{pernu}). Therefore, the functions $J_s$
are periodic functions of $\nu$. Hence,
\beq\label{perJF}
F_i(z, \nu+N)=F_i(z,\nu).
\eeq
For each $\nu$ the space of functions spanned by
the abelian functions $F_i(z,\nu)$ is finite-dimensional. Due to periodicity
of $F_i$ in $\nu$ the total space $\F$ spanned by sequences $F_i(z,\nu)$
is also finite-dimensional. Let $\{F_{\a}\mid\a\in A\}$, for finite set $A$, be
a basis of the factor- space of $\F$ modulo $z$-independent sequences. Then for
all $i\notin A$ there exist constants $c_{i,\a}, d_i(\nu)$ such that
\beq\label{sh03}
F_{i}(z,\nu)-\sum_{\a\in A} c_{i,\a}F_{\a}(z,\nu)=d_{i}(\nu)\,.
\eeq
The rest of the proof of Theorem 1.2 is identical to that in the proof of Theorem 1.1.
Namely,
\begin{lem} Let $\psi$ be a wave function corresponding to $u$, and let
$L_i,\ i\notin A$  be the difference
operator given by the formula
\beq\label{a2d}
L_i=\L^i_+-\sum_{\a\in A} c_{i,\a}\L^{\a}_+, \ i\notin A,
\eeq
where the constants $c_{i,\a}$ are defined by equations (\ref{sh03}).

Then the equations
\beq\label{lpd}
L_i\,\psi=a_i(k)\,\psi, \ \ \ a_i(k)=k^i+\sum_{s=1}^{\infty}a_{s,i}k^{n-s}\,,
\eeq
where $a_{s,i}$ are constants, hold.
\end{lem}
{\it Proof.} From (\ref{m25}) it follows that
\beq\label{lax3d}
[T_1-T-u,L_i]=0.
\eeq
Hence, if $\psi$ is the  wave solution of (\ref{laxm1})
then $L_i\psi$ is also a wave solution
of the same equation. By uniqueness of the wave function up to a constant in $z$-factor
we get (\ref{lp}) and thus the lemma is proven.

\begin{cor} The operators $L_i^{z}$
commute with each other,
\beq\label{com1d}
[L_i^{z},L_j^{z}]=0\,.
\eeq
\end{cor}



\begin{thebibliography}{9}
\bibitem{amkm} H. Airault, H. McKean, and J. Moser.
\emph{Rational and elliptic solutions of the Korteweg-de~Vries
equation and related many-body problem}, Commun.\ Pure Appl.\
Math., {\bf 30} (1977), no.\,1, 95--148.

\bibitem{arb-decon}
E. Arbarello, C. De Concini, {\it On a set of equations characterizing Riemann matrices.}
Ann. of Math. (2) 120 (1984), no. 1, 119--140.

\bibitem{arbarello}
E. Arbarello, C. De Concini, {\it
Another proof of a conjecture of S.P. Novikov on periods of
abelian integrals on Riemann surfaces}, Duke Math. Journal, 54 (1987), 163--178.

\bibitem{krbab} O. Babelon, E. Billey, I. Krichever and M. Talon.
\emph{Spin generalisation of the Calogero-Moser system and
the matrix KP equation}, in ``Topics in Topology and Mathematical
Physics'', Amer. Math. Soc. Transl. Ser.\,2 {\bfseries 170}, Amer.
Math. Soc., Providence, 1995, 83--119.

\bibitem{benzvi}
David Ben-Zvi, Thomas Nevins,  Flows of Calogero-Moser Systems.
arXiv:math/0603722

\bibitem{ch1}
J.L.Burchnall, T.W. Chaundy, {\it Commutative ordinary differential
operators.I}, Proc. London Math Soc. {\bf 21} (1922), 420--440.

\bibitem{ch2}
J.L.Burchnall, T.W. Chaundy, {\it Commutative ordinary differential
operators.II}, Proc. Royal Soc. London  {\bf 118} (1928), 557--583.

\bibitem{fay}
J.D. Fay, {\it Theta functions on Riemann surfaces.}
Lecture Notes in Mathematics, Vol. 352. Springer-Verlag, Berlin-New York, 1973.

\bibitem{gun1}
R. Gunning, {\it Some curves in abelian varieties,} Invent. Math. 66 (1982), no. 3,
377--389.

\bibitem{nekr}
Gorsky, A.; Nekrasov, N. Hamiltonian systems of Calogero-type, and
two-dimensional Yang-Mills theory. Nuclear Phys. B 414 (1994), no.
1--2, 213--238.

\bibitem{kr1}
I.~M. Krichever, {\it Integration of non-linear equations by methods of algebraic geometry},
Funct. Anal. Appl., 11 (1977), n. 1, 12--26.

\bibitem{kr2}
I.~M. Krichever, {\it Methods of algebraic geometry in the theory of non-linear equations},
Russian Math. Surveys, 32 (1977), n. 6, 185--213.

\bibitem{kr-dif}
I.Krichever, {\it Algebraic curves and non-linear difference equation},
Uspekhi Mat. Nauk 33 (1978), n 4, 215-216.

\bibitem{krelkp} I. Krichever.
\emph{Elliptic solutions of Kadomtsev-Petviashvili equations
and integrable systems of particles} (In Russian), Funct. Anal.
Appl., {\bfseries 14} (1980), no.\,1, 45--54. Transl.
282--290.

\bibitem{kr-toda}
I.Krichever, {\it  The periodic nonabelian Toda lattice
and two-dimensional generalization}\,,
appendix to: B.Dubrovin,{\it Theta-functions and nonlinear equations}\,,
Uspekhi Mat. Nauk {\bf 36} (1981), n 2,
72-77.

\bibitem{krdd}
I. Krichever, {\it Two-dimensional periodic difference operators and algebraic
geometry}, Doklady Akad. Nauk USSR 285 (1985), n 1, 31-36.

\bibitem{krnest} I. Krichever.
\emph{Elliptic solutions to difference non-linear equations and
nested Bethe ansatz equations}, Calogero-Moser-Sutherland models
(Montreal, QC, 1997), 249--271, CRM Ser. Math. Phys., Springer,
New York, 2000.

\bibitem{kreltoda} I. Krichever.
\emph{Elliptic analog of the Toda lattice}, Internat. Math. Res.
Notices (2000), no. 8, 383--412.

\bibitem{krlax}
Krichever, I. Vector bundles and Lax equations on algebraic curves.
Comm. Math. Phys. 229 (2002), no. 2, 229--269.

\bibitem{kr-schot}
I.Krichever, {\it Integrable linear equations and the
Riemann-Schottky problem},  Algebraic geometry and number theory,
497--514, Progr. Math., 253, Birkh\"auser Boston, Boston, MA, 2006.

\bibitem{kr-tri}
I. Krichever, {\it Characterizing Jacobians via trisecants of the
Kummer Variety}, math.AG /0605625.

\bibitem{krwz} I. Krichever, O. Lipan, P. Wiegmann, and A. Zabrodin.
\emph{Quantum integrable models and discrete classical Hirota
equations}, Comm. Math. Phys., {\bfseries 188} (1997), no.\,2,
267--304.


\bibitem{KrSi}
I.Krichever, T.Shiota, {\it Abelian solutions of the KP equation}, arXiv:0804.0274

\bibitem{krzab} I. Krichever, A. Zabrodin.
\emph{Spin generalisation of the Ruijsenaars-Schneider model,
the nonabelian two-dimensionalized Toda lattice, and
representations of the Sklyanin algebra} (in Russian), Uspekhi
Mat. Nauk, {\bfseries 50} (1995), no.\,6, 3--56.

\bibitem{shiota}
T. Shiota, {\it Characterization of {Jacobian} varieties in terms of soliton
equations}, Invent. Math., 83(2):333--382, 1986.

\bibitem{mum}
D. Mumford, {\it An algebro-geometric construction of commuting operators and of solutions to
the Toda lattice  equation, Korteweg-de Vries equation and related non-linear equations} --
Proceedings Int.Symp. Algebraic Geometry, Kyoto, 1977, 115--153, Kinokuniya Book Store,
Tokyo, 1978.

\bibitem{ruij}
S.N.M. Ruijsenaars and H. Schneider,
\emph{A new class of integrable systems and its relation to solitons},
Ann.\ Physics 170 (1986) 370--405.

\bibitem{serr}
J-P. Serre, {\it Faisceaux alg\'ebriques coh\'erents}, (French) Ann. of Math. (2)
61, (1955). 197--278.

\bibitem{wel}
G.E. Welters, {\it On flexes of the Kummer variety (note on a theorem of R. C. Gunning).}
Nederl. Akad. Wetensch. Indag. Math. 45 (1983), no. 4, 501--520.

\bibitem{wel1}
G.E. Welters, {\it A criterion for {Jacobi} varieties}, Ann. of Math.,
120 (1984), n. 3, 497--504.

\end{thebibliography}
\end{document}